\documentclass{article}
\usepackage{amsmath}
\usepackage{amsfonts}

\input{tcilatex}
\begin{document}

\begin{center}
\textbf{LINEAR WEINGARTEN FACTORABLE SURFACES IN ISOTROPIC SPACES}

\bigskip

Muhittin Evren Aydin, Alper Osman Ogrenmis

Firat University, Turkey

\bigskip
\end{center}

\textbf{Abstract. }In this paper, we deal with the linear Weingarten
factorable surfaces in the isotropic 3-space $\mathbb{I}^{3}$ satisfying the
relation $aK+bH=c,$ where $K$ is the relative curvature and $H$ the
isotropic mean curvature, $a,b,c\in 
\mathbb{R}
$. We obtain a complete classification for such surfaces in $\mathbb{I}^{3}.$
As a further study, we classify all graph surfaces in $\mathbb{I}^{3}$
satisfying the relation $K=H^{2},$ which is the equality case of the famous
Euler inequality for surfaces in a Euclidean space.

\bigskip



\section{Introduction}

Let $M^{2}$ be a regular surface of a Euclidean 3-space $\mathbb{%
\mathbb{R}
}^{3}$ and $\kappa _{1},\kappa _{2}$ principal curvatures of $M^{2}.$\ A
surface $M^{2}$ in $%
\mathbb{R}
^{3}$ is called a \textit{Weingarten surface (W-surface)} if it satisfies
the following non-trivial functional relation 
\begin{equation}
\phi \left( \kappa _{1},\kappa _{2}\right) =0  \tag{1.1}
\end{equation}%
for a smooth function $\phi $ of two variables. $\left( 1.1\right) $
immediately yields%
\begin{equation}
\delta \left( K,H\right) =0,  \tag{1.2}
\end{equation}%
where $K$ and $H$ are respectively the Gaussian and mean curvatures of $%
M^{2} $.

$\left( 1.2\right) $ is equivalent to the vanishing of the corresponding
Jacobian determinant, i.e. $\left\vert \partial \left( K,H\right) /\partial
\left( u,v\right) \right\vert =0$ for a coordinate pair $\left( u,v\right) $
on $M^{2}.$

If $M^{2}$ is a surface in $%
\mathbb{R}
^{3}$ verifying the following relation 
\begin{equation}
aH+bK=c,\text{ }a,b,c\in 
\mathbb{R}
,\text{ }\left( a,b,c\right) \neq \left( 0,0,0\right) ,  \tag{1.3}
\end{equation}%
then it is called a\textit{\ linear Weingarten surface (LW-surface). }If $a=0
$ or $b=0$ in $\left( 1.3\right) ,$ then the LW-surfaces reduce to the ones
with constant curvature. Many geometers extensively have studied such
surfaces, see \cite{8,10,17}, \cite{19}-\cite{21}, \cite{37}.

On the other hand, let $M^{2}$ be a graph surface of a smooth function $%
z=z\left( x,y\right) $. If $z\left( x,y\right) =f\left( x\right) g\left(
y\right) ,$ then $M^{2}$ is called a \textit{factorable surface} or \textit{%
homothetical surface}. For geometric results on these surfaces in ambient
spaces, see \cite{1}-\cite{4},\cite{11,14,22,35,36}.

Most recently, the first author and M. Ergut \cite{1} classified the
factorable surfaces with constant relative and constant isotropic mean
curvature in the isotropic 3-space $\mathbb{I}^{3}$ which has been
introduced by K. Strubecker \cite{34} and H. Sachs \cite{31,32}.

Several classes of surfaces in $\mathbb{I}^{3}$ have been studied by I.
Kamenarovic (\cite{15}), B. Pavkovic (\cite{28}), Z. M. Sipus (\cite{33})
and M.K. Karacan and et al. (\cite{16}).

The main goal of this paper is to study LW-factorable surfaces in $\mathbb{I}%
^{3}.$ In the present paper, we provide a classification for the
LW-factorable surfaces in $\mathbb{I}^{3}$. As a further study, we classify
the graph surfaces in $\mathbb{I}^{3}$ satisfying the relation $K=H^{2}.$

\section{Preliminaries}

The isotropic 3-space $\mathbb{I}^{3}$ is a Cayley--Klein space defined from
a 3-dimensional projective space $P\left( 
\mathbb{R}
^{3}\right) $ with the absolute figure which is an ordered triple $\left(
\omega ,f_{1},f_{2}\right) $, where $\omega $ is a plane in $P\left( 
\mathbb{R}
^{3}\right) $ and $f_{1},f_{2}$ are two complex-conjugate straight lines in $%
\omega $. For more details, we refer \cite{5,7,9,27,29,30,32}

The homogeneous coordinates in $P\left( 
\mathbb{R}
^{3}\right) $ are introduced in such a way that the \textit{absolute plane} $%
\omega $ is given by $X_{0}=0$ and the \textit{absolute lines} $f_{1},f_{2}$
by $X_{0}=X_{1}+iX_{2}=0,$ $X_{0}=X_{1}-iX_{2}=0.$ The intersection point $%
F(0:0:0:1)$ of these two lines is called the \textit{absolute point}. Affine
coordinates in $P\left( 
\mathbb{R}
^{3}\right) $ are given by 
\begin{equation*}
x_{1}=\frac{X_{1}}{X_{0}},\text{ }x_{2}=\frac{X_{2}}{X_{0}},\text{ }x_{3}=%
\frac{X_{3}}{X_{0}}.
\end{equation*}

Consider the points $x=\left( x_{1},x_{2},x_{3}\right) $ and $y=\left(
y_{1},y_{2},y_{3}\right) .$ \textit{Isotropic distance }$d_{\mathbb{I}%
}\left( x,y\right) $\textit{\ }of two points $x$ and $y$ is defined as%
\begin{equation*}
d_{\mathbb{I}}\left( x,y\right) =\left( y_{1}-x_{1}\right) ^{2}+\left(
y_{2}-x_{2}\right) ^{2}.
\end{equation*}%
The lines in $x_{3}-$direction are called \textit{isotropic} lines. The
plane containing an isotropic line is called an\textit{\ isotropic plane. }%
Other planes\textit{\ }are\textit{\ non-isotropic.}

Let $M^{2}$ be a graph surface in $\mathbb{I}^{3}$ corresponding to a smooth
function $z=z\left( x,y\right) $ on a open domain $D\subseteq \mathbb{R}^{2}$%
. Then it is parameterized as follows

\begin{equation}
r:D\subseteq \mathbb{R}^{2}\longrightarrow \mathbb{I}^{3}:\text{ }\left(
x,y\right) \longmapsto \left( x,y,z\left( x,y\right) \right) .  \tag{2.1}
\end{equation}%
It follows from $\left( 2.1\right) $ that $M^{2}$ is an admissble surface
immersed in $\mathbb{I}^{3}$ (i.e. without isotropic tangent planes). The
reader can find a well bibliography for surfaces of $\mathbb{I}^{3}$ in \cite%
{32}.

The metric on $M^{2}$ induced from $\mathbb{I}^{3}$ is given by $g_{\ast
}=dx^{2}+dy^{2}.$ This implies that $M^{2}$ is always flat with respect to
the induced metric $g_{\ast }.$ Thus its Laplacian is given by%
\begin{equation*}
\bigtriangleup =\frac{\partial ^{2}}{\partial x^{2}}+\frac{\partial ^{2}}{%
\partial y^{2}}.
\end{equation*}%
The \textit{relative curvature} $K$ and the \textit{isotropic mean curvature 
}$H$\textit{\ of }$M^{2}$ are respectively defined by%
\begin{equation}
K=z_{xx}z_{yy}-\left( z_{xy}\right) ^{2}  \tag{2.2}
\end{equation}%
and%
\begin{equation}
H=\bigtriangleup z=\frac{z_{xx}+z_{yy}}{2}.  \tag{2.3}
\end{equation}%
A surface is called \textit{isotropic minimal }(resp. \textit{isotropic flat}%
) if $H$ (resp. $K$) vanishes.

\section{LW-factorable surfaces in $\mathbb{I}^{3}$}

Let $M^{2}$ be a factorable surface in $\mathbb{I}^{3}.$ Then it is a graph
surface of a smooth function $z\left( x,y\right) =f\left( x\right) g\left(
y\right) .$ By $\left( 2.2\right) $ and $\left( 2.3\right) ,$ the relative
and isotropic mean curvatures of $M^{2}$ respectively turn to 
\begin{equation}
K=\left( f^{\prime \prime }f\right) \left( g^{\prime \prime }g\right)
-\left( f^{\prime }\right) ^{2}\left( g^{\prime }\right) ^{2}  \tag{3.1}
\end{equation}%
and 
\begin{equation}
2H=f^{\prime \prime }g+fg^{\prime \prime },  \tag{3.2}
\end{equation}%
where $f^{\prime }=\frac{df}{dx}$ and $g^{\prime }=\frac{dg}{dy},$ etc.

We mainly aim to classify the LW-factorable surfaces in $\mathbb{I}^{3}$.
For this, let $M^{2}$ be a LW-factorable surface in $\mathbb{I}^{3}$
satisfying the relation $\left( 1.3\right) $. Since at least one of $a,$ $b$
and $c$ is nonzero in $\left( 1.3\right) $, without loss of generality, we
may assume $b\neq 0.$ By dividing both sides of $\left( 1.3\right) $ with $b$
and putting $\frac{a}{b}=2m_{0}$ and $\frac{c}{b}=n_{0},$ we write%
\begin{equation}
2m_{0}H+K=n_{0},\text{ }m_{0},n_{0}\in 
\mathbb{R}
.  \tag{3.3}
\end{equation}%
When $m_{0}=0$ in $\left( 3.3\right) ,$ $M^{2}$ becomes a factorable surface
in $\mathbb{I}^{3}$ with $K=const.,$ however, such surfaces were already
classified in \cite{1}. In our framework, it is meaningful to take $%
m_{0}\neq 0.$

By $\left( 3.1\right) -\left( 3.3\right) ,$ we get%
\begin{equation}
\left( f^{\prime \prime }f\right) \left( g^{\prime \prime }g\right) -\left(
f^{\prime }\right) ^{2}\left( g^{\prime }\right) ^{2}+m_{0}\left( f^{\prime
\prime }g+fg^{\prime \prime }\right) =n_{0}.  \tag{3.4}
\end{equation}%
We have to distinguish some situations in order to solve $\left( 3.4\right)
. $

\bigskip

\textbf{Remark 3.1. }\textit{From now on, we use the notation }$c_{i}$%
\textit{\ to denote nonzero constants and }$d_{i}$\textit{\ to denote some
constants, }$i=1,2,3,...$

\bigskip

\textbf{Case 1. }$f\left( x\right) =f_{0}\in 
\mathbb{R}
-\left\{ 0\right\} .$ By $\left( 3.4\right) $, we find%
\begin{equation}
g\left( y\right) =\frac{n_{0}}{f_{0}m_{0}}y^{2}+d_{1}y+d_{2}.  \tag{3.5}
\end{equation}%
Similarly, it can be obtained from $\left( 3.4\right) $ that $f\left(
x\right) =\frac{n_{0}}{g_{0}m_{0}}x^{2}+d_{3}x+d_{4}$ when $g\left( y\right)
=g_{0}\in 
\mathbb{R}
-\left\{ 0\right\} .$

\bigskip

\textbf{Remark 3.2. }\textit{In Case 1 (i.e. in the case }$f\left( x\right)
=f_{0}$ or $g\left( y\right) =g_{0},$ $f_{0},g_{0}\in 
\mathbb{R}
-\left\{ 0\right\} $)\textit{, }$M^{2}$ \textit{is an isotropic flat
factorable surface} \textit{in} $\mathbb{I}^{3}$ \textit{with} $H=\frac{n_{0}%
}{m_{0}}.$

\bigskip

\textbf{Case 2. }Let $f$ be a linear function, i.e. $f\left( x\right)
=c_{1}x+d_{5}.$ It follows from $\left( 3.4\right) $ that%
\begin{equation}
-c_{1}^{2}\left( g^{\prime }\right) ^{2}+m_{0}\left\{ \left(
c_{1}x+d_{5}\right) g^{\prime \prime }\right\} =n_{0}.  \tag{3.6}
\end{equation}%
Taking partial derivative of $\left( 3.6\right) $ with respect to $x$ gives $%
m_{0}c_{1}g^{\prime \prime }=0,$ namely $g\left( y\right) =c_{2}y+d_{6}.$
With similar arguments, we can find if $g$ is a linear function in $\left(
3.4\right) $, so is $f.$

\bigskip

\textbf{Remark 3.3. }\textit{In Case 2 (i.e. in the case }$f\left( x\right)
=c_{1}x+d_{5}$ and $g\left( y\right) =c_{2}y+d_{6}$\textit{), }$M^{2}$ 
\textit{is an isotropic minimal factoable surface in} $\mathbb{I}^{3}$ 
\textit{with }$K=-\left( c_{1}c_{2}\right) ^{2}$.

\bigskip

\textbf{Case 3. }$f$ and $g$ are non-linear functions. By dividing $\left(
3.4\right) $ with the product $ff^{\prime \prime },$ we have%
\begin{equation}
g^{\prime \prime }g-\frac{\left( f^{\prime }\right) ^{2}}{ff^{\prime \prime }%
}\left( g^{\prime }\right) ^{2}+m_{0}\frac{g}{f}+m_{0}\frac{g^{\prime \prime
}}{f^{\prime \prime }}=\frac{n_{0}}{ff^{\prime \prime }}.  \tag{3.7}
\end{equation}%
By taking partial derivative $\left( 3.7\right) $ with respect to $y$ and
then dividing with $g^{\prime }g^{\prime \prime }$, we deduce%
\begin{equation}
1+\frac{gg^{\prime \prime \prime }}{g^{\prime }g^{\prime \prime }}-2\frac{%
\left( f^{\prime }\right) ^{2}}{ff^{\prime \prime }}+\left( \frac{m_{0}}{f}%
\right) \frac{1}{g^{\prime \prime }}+\left( \frac{m_{0}}{f^{\prime \prime }}%
\right) \frac{g^{\prime \prime \prime }}{g^{\prime }g^{\prime \prime }}=0. 
\tag{3.8}
\end{equation}%
We have two cases:

\bigskip

\textbf{Case 3.1. }$g^{\prime \prime \prime }=0,$ i.e. 
\begin{equation}
g\left( y\right) =c_{3}y^{2}+d_{7}y+d_{8}.  \tag{3.9}
\end{equation}%
Then $\left( 3.8\right) $ reduces to%
\begin{equation}
1-2\frac{\left( f^{\prime }\right) ^{2}}{ff^{\prime \prime }}+\left( \frac{%
m_{0}}{2c_{3}}\right) \frac{1}{f}=0.  \tag{3.10}
\end{equation}%
$\left( 3.10\right) $ can be rewritten as%
\begin{equation}
\left( \frac{m_{0}}{2c_{3}}+f\right) f^{\prime \prime }-2\left( f^{\prime
}\right) ^{2}=0.  \tag{3.11}
\end{equation}%
After solving $\left( 3.11\right) ,$ we find%
\begin{equation}
f\left( x\right) =-\left( \frac{1}{c_{4}x+d_{9}}+\frac{m_{0}}{2c_{3}}\right)
.  \tag{3.12}
\end{equation}%
Considering $\left( 3.9\right) $ and $\left( 3.12\right) $ into $\left(
3.4\right) $ gives that%
\begin{equation}
\frac{c_{4}^{2}}{\left( c_{4}x+d_{9}\right) ^{4}}\left(
4c_{3}d_{8}-d_{7}^{2}\right) -\frac{2m_{0}c_{3}}{c_{4}x+d_{9}}%
-m_{0}^{2}=n_{0}.  \tag{3.13}
\end{equation}%
In the particular case $d_{7}=d_{8}=0,$ we obtain the following contradiction%
\begin{equation*}
x=-\frac{1}{c_{4}}\left( \frac{2m_{0}c_{3}}{n_{0}+m_{0^{2}}}+d_{9}\right)
\end{equation*}%
since $x$ is an independent variable.

\bigskip

\textbf{Case 3.2. }$g^{\prime \prime \prime }\neq 0.$ By taking partial
derivatives of $\left( 3.8\right) $ with respect to $x$ and $y$, we conclude 
\begin{equation}
\left( \frac{f^{\prime }}{f^{2}}\right) \frac{g^{\prime \prime \prime }}{%
\left( g^{\prime \prime }\right) ^{2}}-\frac{f^{\prime \prime \prime }}{%
\left( f^{\prime \prime }\right) ^{2}}\left( \frac{g^{\prime \prime \prime }%
}{g^{\prime }g^{\prime \prime }}\right) ^{\prime }=0.  \tag{3.14}
\end{equation}%
Since $f^{\prime }\neq 0\neq g^{\prime \prime \prime }\ $, neither $%
f^{\prime \prime \prime }$ nor $\left( \frac{g^{\prime \prime \prime }}{%
g^{\prime }g^{\prime \prime }}\right) ^{\prime }$ can vanish in $\left(
3.14\right) $. Then $\left( 3.14\right) $ can be rewritten as%
\begin{equation}
\frac{f^{\prime }\left( f^{\prime \prime }\right) ^{2}}{f^{2}f^{\prime
\prime \prime }}=\frac{\left( g^{\prime \prime }\right) ^{2}}{g^{\prime
\prime \prime }}\left( \frac{g^{\prime \prime \prime }}{g^{\prime }g^{\prime
\prime }}\right) ^{\prime }.  \tag{3.15}
\end{equation}%
Since the left side of $\left( 3.15\right) $ is a function of $x,$ however
the right side is a function of $y.$ Then both sides have to be equal a
nonzero constant, i.e.%
\begin{equation}
\frac{f^{\prime }\left( f^{\prime \prime }\right) ^{2}}{f^{2}f^{\prime
\prime \prime }}=c_{5}=\frac{\left( g^{\prime \prime }\right) ^{2}}{%
g^{\prime \prime \prime }}\left( \frac{g^{\prime \prime \prime }}{g^{\prime
}g^{\prime \prime }}\right) ^{\prime }.  \tag{3.16}
\end{equation}%
\newline
From the left side of $\left( 3.16\right) ,$ we write%
\begin{equation}
\frac{f^{\prime \prime \prime }}{\left( f^{\prime \prime }\right) ^{2}}=%
\frac{1}{c_{5}}\frac{f^{\prime }}{f^{2}}  \tag{3.17}
\end{equation}%
or, by taking once integral with respect to $x$,%
\begin{equation}
f^{\prime \prime }=\frac{c_{5}f}{c_{5}d_{10}f+1}  \tag{3.18}
\end{equation}%
for an integration constant $d_{10}.$ Assuming $d_{10}=0$ in $\left(
3.18\right) $ gives $f^{\prime \prime }=c_{5}f.$ By putting this in $\left(
3.4\right) $ we derive%
\begin{equation}
\left( f^{2}\right) \left( g^{\prime \prime }g\right) -\left( f^{\prime
}\right) ^{2}\left( g^{\prime }\right) ^{2}+m_{0}f\left( c_{5}g+g^{\prime
\prime }\right) =n_{0}.  \tag{3.19}
\end{equation}%
Dividing $\left( 3.19\right) $ with $f$ and then taking partial derivative
with respect to $x$ imply%
\begin{equation}
g^{\prime \prime }g-\left\{ 2\frac{f^{\prime \prime }}{f}-\left( \frac{%
f^{\prime }}{f}\right) ^{2}\right\} \left( g^{\prime }\right) ^{2}=\frac{%
-n_{0}}{f^{2}}.  \tag{3.20}
\end{equation}%
If $2\frac{f^{\prime \prime }}{f}-\left( \frac{f^{\prime }}{f}\right) ^{2}$
is some constant in $\left( 3.20\right) ,$ then, by taking a partial
derivative of $\left( 3.20\right) $ with respect to $x$, we obtain%
\begin{equation*}
0=\frac{2n_{0}f^{\prime }}{f^{3}},
\end{equation*}%
which is not possible since $f$ is non-linear. Now by again taking partial
derivative of $\left( 3.20\right) $ with respect to $x,$ we deduce%
\begin{equation}
-\left\{ 2\frac{f^{\prime \prime }}{f}-\left( \frac{f^{\prime }}{f}\right)
^{2}\right\} ^{\prime }\frac{f^{3}}{2n_{0}f^{\prime }}=\frac{1}{\left(
g^{\prime }\right) ^{2}}.  \tag{3.21}
\end{equation}%
Since $g^{\prime \prime \prime }\neq 0$, the right side of $\left(
3.21\right) $ is a function of $y,$ but the left side is either a nonzero
constant or a function of $x$. Both cases are not possible.

Therefore we have proved the following:

\bigskip

\textbf{Theorem 3.1. }\textit{Let }$M^{2}$ \textit{be a LW-factorable
surface which is the graph of }$z\left( x,y\right) =f\left( x\right) g\left(
y\right) $ \textit{in }$\mathbb{I}^{3}$\textit{. Then we have one of the
following statements:}

\textit{(A) }$f\left( x\right) =f_{0}\in 
\mathbb{R}
-\left\{ 0\right\} ,$ $g\left( y\right) =c_{6}y^{2}+d_{11}y+d_{12};$

\textit{(B) }$g\left( y\right) =g_{0}\in 
\mathbb{R}
-\left\{ 0\right\} ,$ $f\left( x\right) =c_{7}x^{2}+d_{13}x+d_{14};$

\textit{(C) }$z\left( x,y\right) =\left( c_{8}x+d_{15}\right) \left(
c_{9}y+d_{16}\right) .$

\bigskip

By Remark 3.2, Remark 3.3 and Theorem 3.1, we immediately derive the
following.

\bigskip

\textbf{Corollary 3.1. }\textit{The LW-factorable surfaces in }$\mathbb{I}%
^{3}$\textit{\ are only the ones whose both }$K$ and $H$ are constants$.$

\section{Graph surfaces with $K=H^{2}$}

Let $M^{2}$ be a surface of the Euclidean 3-space $\mathbb{R}^{3}$. The
Euler inequality for $M^{2}$ including the Gaussian $\left( K\right) $ and
mean curvature $\left( H\right) $ follows%
\begin{equation}
K\leq H^{2}.  \tag{4.1}
\end{equation}%
For more generalizations of this inequality, see \cite{6}, \cite{24}-\cite%
{26}.

The equality sign of $\left( 4.1\right) $ holds on $M^{2}$ if and only if it
is totally umbilical, i.e. a part of a plane or a two sphere in $\mathbb{E}%
^{3}$.

Now we are interested with the factorable surfaces in $\mathbb{I}^{3}$
satisfying $K=H^{2}.$ For this aim, let us reconsider $\left( 3.1\right) $
and $\left( 3.2\right) $. If $K=H^{2},$ then%
\begin{equation}
\left( f^{\prime \prime }g-fg^{\prime \prime }\right) ^{2}+4\left( f^{\prime
}g^{\prime }\right) ^{2}=0.  \tag{4.2}
\end{equation}%
$\left( 4.2\right) $ immediately implies that%
\begin{equation}
f^{\prime \prime }g-fg^{\prime \prime }=0\text{ and }f^{\prime }g^{\prime
}=0.  \tag{4.3}
\end{equation}%
By $\left( 4.3\right) $ we conclude that either $f=const$. and $g\left(
y\right) =c_{1}y+d_{1}$ or $g=const$. and $f\left( x\right) =c_{2}x+d_{2}.$
It yields that a factorable surface satisfying $K=H^{2}$ is a non-isotropic
plane in $\mathbb{I}^{3}$.

Therefore we have proved the following:

\bigskip

\textbf{Proposition 4.1. }\textit{The factorable surfaces in} $\mathbb{I}%
^{3} $ \textit{satisfying} $K=H^{2}$ \textit{are only non-isotropic planes.}

\bigskip

As a generalization of Proposition 4.1, we are able to investigate all graph
surfaces in $\mathbb{I}^{3}$ satisfying $K=H^{2}.$ More precisely, let $%
M^{2} $ be a graph surface of the smooth function $z=z\left( x,y\right) $ in 
$\mathbb{I}^{3}.$ If $K=H^{2}$ on $M^{2},$ then we get%
\begin{equation}
\left( z_{xx}-z_{yy}\right) ^{2}+4\left( z_{xy}\right) ^{2}=0.  \tag{4.4}
\end{equation}%
$\left( 4.4\right) $ yields that%
\begin{equation}
z_{xy}=0  \tag{4.5}
\end{equation}%
and%
\begin{equation}
z_{xx}=z_{yy}.  \tag{4.6}
\end{equation}%
By $\left( 4.5\right) ,$ we derive%
\begin{equation}
z\left( x,y\right) =\alpha \left( x\right) +\beta \left( y\right)  \tag{4.7}
\end{equation}%
and considering $\left( 4.7\right) $ into $\left( 4.6\right) $ gives%
\begin{equation}
\frac{d^{2}\alpha }{dx^{2}}=\frac{d^{2}\beta }{dy^{2}}=d_{3},\text{ }%
d_{3}\in 
\mathbb{R}
.  \tag{4.8}
\end{equation}%
By solving $\left( 4.8\right) ,$ we find%
\begin{equation}
\alpha \left( x\right) =\frac{d_{3}}{2}x^{2}+d_{4}x+d_{5},\text{ }\beta
\left( y\right) =\frac{d_{3}}{2}y^{2}+d_{6}y+d_{7}.  \tag{4.9}
\end{equation}%
$\left( 4.9\right) $ implies that $M^{2}$ is either a non-isotropic plane $%
\left( d_{3}=0\right) $ or a parabolic sphere $\left( d_{3}\neq 0\right) $
of $\mathbb{I}^{3}.$ For more details of planes and spheres in $\mathbb{I}%
^{3}$, see \cite{29,32}.

Consequently, we have

\bigskip

\textbf{Theorem 4.1. }\textit{A graph surface of a function }$z=z\left(
x,y\right) $ \textit{in} $\mathbb{I}^{3}$ \textit{with} $K=H^{2}$ \textit{is
either (a piece of) a non-isotropic plane or (a piece of) a parabolic sphere
given by}%
\begin{equation*}
z\left( x,y\right) =c_{3}\left( x^{2}+y^{2}\right) +d_{8}x+d_{9}y+d_{10}.
\end{equation*}

\bigskip

M.E. Aydin

Department of Mathematics

Firat University

23119 Elazig

Turkey

E-mail: meaydin@firat.edu.tr

\bigskip

A. O.Ogrenmis

Department of Mathematics

Firat University

23119 Elazig

Turkey

E-mail: aogrenmis@firat.edu.tr

\end{document}